# Recurrences for some sequences of binomial sums II: A simpler approach

*Johann Cigler*

Fakultät für Mathematik, Universität Wien

Due to some results by John P. D'Angelo and Dusty Grundmeier about CR- mappings the main results of my 2001 paper about recurrences for some sequences of binomial sums can be simplified.

## 1. Introduction

The main result of my paper [1] is

**Lemma 1**
*For each pair of integers* $n \geq 1, m \geq 1$ *there exist uniquely determined integers* $a(n,m,j)$ *such that*

$$(x+1)^n + \sum_{k=1}^{m-1} \sum_{j=1+\left\lfloor \frac{(k-1)n}{m} \right\rfloor}^{\left\lfloor \frac{kn}{m} \right\rfloor} a(n,m,j)(x+1)^{kn-mj} x^j = 1 + (-1)^{m(n-1)} x^n. \tag{1.1}$$

To formulate (1.1) in a more succinct way define polynomials $p_k(n,m,x,s)$ by

$$p_1(n,m,x,s) = x^n + \sum_{j=1}^{\left\lfloor \frac{n}{m} \right\rfloor} a(n,m,j) x^{n-mj} s^j \tag{1.2}$$

and

$$p_k(n,m,x,s) = \sum_{j=1+\left\lfloor \frac{(k-1)n}{m} \right\rfloor}^{\left\lfloor \frac{kn}{m} \right\rfloor} a(n,m,j) x^{kn-mj} s^j \tag{1.3}$$

for $n \geq 1$.

Then (1.1) becomes for $m \geq 2$

$$\sum_{k=1}^{m-1} p_k(n,m,x+1,x) = 1 + (-1)^{m(n-1)} x^n. \tag{1.4}$$

**Attention:** To make some formulae more elegant we define $p_k(n,m,x,s)$ also for $n=0$ as $p_k(0,m,x,s) = (-1)^{k-1}\binom{m}{k}$. Then e.g. (1.4) remains true for $n=0$. But in some cases we need that $p_1(n,m,x,s)$ takes the value $1$ for $n=0$. To this end we define polynomials $\tilde{p}_1(n,m,x,s)$ by $\tilde{p}_1(n,m,x,s) = p_1(n,m,x,s)$ for $n \geq 1$ and $\tilde{p}_1(0,m,x,s) = 1$.

It is remarkable that the same sort of polynomials occurred independently in the study of certain CR mappings in the theory of several complex variables.

In 2004 John P. D'Angelo [3] defined polynomials $f_{n,m}(x,s)$, which will be called "invariant polynomials", by

$$f_{n,m}(x,s)=1-\prod_{j=0}^{n-1}\left(1-\omega^j x+\omega^{mj}s\right) \tag{1.5}$$

where $\omega$ denotes a primitive $n-$th root of unity. (In fact he used a slightly different notation).
He showed that $f_{n,m}(x,s)$ is the unique polynomial $f(x,s)$ satisfying

1) $f(0,0)=0$
2) $f(x,s)=1$ when $x-s=1$
3) $\deg f=n$
4) the invariance property $f(\omega x,\omega^m s)=f(x,s)$ for all $x$ and $s$.

A polynomial $m(x,s)$ has weight $k$ if $m(\lambda x,\lambda^m s)=\lambda^{kn}m(x,s)$.

This characterization shows that

$$f_{n,m}(x,s)=\sum_{k=1}^{m}p_k(n,m,x,s). \tag{1.6}$$

For 1) and 3) are obvious from the definition. Condition 2) is (1.4) if we define $p_m(n,m,x,s)=-(-1)^{(m-1)n}s^n$.
The polynomials $p_k(n,m,x,s)$ have weight $k$ for each $k$ and therefore 4) is also satisfied.

In this note I want to give a simpler approach to my paper [1] taking into account these results.

## 2. A direct approach for small $m$

Let us first make some general observations about (1.1).
Let us assume that there exists a solution $\left(a(n,m,j)\right)_{j=1}^{n}$ of the equation

$$(x+1)^n+\sum_{k=1}^{m}\sum_{j=1+\left\lfloor\frac{(k-1)n}{m}\right\rfloor}^{\left\lfloor\frac{kn}{m}\right\rfloor}a(n,m,j)(x+1)^{kn-mj}x^j-1=0.$$

Then all coefficients of the left-hand side must vanish. This is equivalent with

$$\binom{n}{k}+\sum_{j=0}^{k}a(n,m,j)\binom{-mj \bmod n}{k-j}=0$$

for $1\le k\le n.$
This means that

$$a(n,m,k)=-\binom{n}{k}-\sum_{j=1}^{k-1}a(n,m,j)\binom{-mj \bmod n}{k-j} \tag{2.1}$$

with initial value

$$a(n,m,1)=-n. \tag{2.2}$$

In particular we see that *all coefficients* $a(n,m,j)$ *are integers.*

**a)** For $m=1$ we get the trivial identity $1=(x+1-x)^n=\sum_{j=0}^{n}\binom{n}{j}(-1)^j(x+1)^{n-j}x^j$.

In this case the table $\left(a(n,1,j)\right)$ consists of binomial coefficients

| | | | | |
|---|---|---|---|---|
| 1 | | | | |
| 1 | -1 | | | |
| 1 | -2 | 1 | | |
| 1 | -3 | 3 | -1 | |
| 1 | -4 | 6 | -4 | 1 |

Here we have

$$p_1(n,1,x,s)=(x-s)^n=\sum_{k=0}^{n}(-1)^k\binom{n}{k}s^kx^{n-k}. \tag{2.3}$$

Note that also

$$\sum_{k=0}^{n}\binom{n}{k}p_1(n-k,1,x,s)s^k=\sum_{k=0}^{n}\binom{n}{k}(x-s)^{n-k}s^k=x^n. \tag{2.4}$$

**b)** For $m=2$ the table $\left(a(n,2,j)\right)$ begins with

| | | | | | | |
|---|---|---|---|---|---|---|
| 1 | | | | | | |
| 1 | -1 | | | | | |
| 1 | -2 | -1 | | | | |
| 1 | -3 | 0 | -1 | | | |
| 1 | -4 | 2 | 0 | -1 | | |
| 1 | -5 | 5 | 0 | 0 | -1 | |
| 1 | -6 | 9 | -2 | 0 | 0 | -1 |

If you are familiar with the sequence of Lucas polynomials $L_n(x,s)$ whose first terms are

$\{2,\ x,\ 2s+x^2,\ 3sx+x^3,\ 2s^2+4sx^2+x^4,\ 5s^2x+5sx^3+x^5,\ 2s^3+9s^2x^2+6sx^4+x^6\}$

you will guess that $a(n,2,j)$ is a coefficient of $L_n(x,-s)$.

The Lucas polynomials can be defined by

$$L_n(x,s)=\left(\frac{x+\sqrt{x^2+4s}}{2}\right)^n+\left(\frac{x-\sqrt{x^2+4s}}{2}\right)^n=\sum_{k=0}^{\left\lfloor\frac{n}{2}\right\rfloor}\binom{n-k}{k}\frac{n}{n-k}s^kx^{n-2k}.$$

They satisfy the recurrence $L_n(x,s)=xL_{n-1}(x,s)+sL_{n-2}(x,s)$ with initial values $L_0(x,s)=2$ and $L_1(x,s)=x$ and the identity

$$L_n(x+y,-xy)=x^n+y^n, \tag{2.5}$$

which immediately implies

$$p_1(n,2,x,s)=L_n(x,-s)=\sum_{k=0}^{\left\lfloor\frac{n}{2}\right\rfloor}(-1)^k\binom{n-k}{k}\frac{n}{n-k}s^kx^{n-2k}. \tag{2.6}$$

and

$$(x+1)^n+\sum_{k=1}^{\left\lfloor\frac{n}{2}\right\rfloor}(-1)^k\binom{n-k}{k}\frac{n}{n-k}(x+1)^{n-2k}x^k=1+x^n$$

for $n\geq 1$.

Let us for further use note the generating function of the Lucas polynomials

$$\sum_{n\geq 0} L_n(x,s)z^n = \frac{2-xz}{1-xz-sz^2}. \qquad (2.7)$$

For the modified polynomials $\tilde{p}_1(n,2,x,s)$ there is a well-known analogue of (2.4)

$$\sum_{k=0}^{\left\lfloor \frac{n}{2} \right\rfloor} \binom{n}{k} \tilde{p}_1(n-2k,2,x,s)s^k = x^n. \qquad (2.8)$$

To prove this let $\alpha = \frac{x+\sqrt{x^2+4s}}{2}$ and $\beta = \frac{x-\sqrt{x^2+4s}}{2}$.

Then

$$\sum_{k=0}^{\left\lfloor \frac{n}{2} \right\rfloor} \binom{n}{k} \tilde{p}_1(n-2k,2,x,s)s^k = \sum_{k=0}^{\left\lfloor \frac{n}{2} \right\rfloor} \binom{n}{k} \left(\alpha^{n-2k}+\beta^{n-2k}\right)(\alpha\beta)^k = \sum_{k=0}^{\left\lfloor \frac{n}{2} \right\rfloor} \binom{n}{k} \left(\alpha^{n-k}\beta^k+\alpha^k\beta^{n-k}\right).$$

If $n$ is odd then the right-hand side is $\left(\alpha+\beta\right)^n = x^n.$ For $n$ even the same is true because the term $\binom{2n}{n}$ occurs only once because of $\tilde{p}_1(0,2,x,s)=1.$

**c)** For $m=3$ the first terms of the sequence $\left(p_1(n,3,x,s)\right)_{n\geq 0}$ are

$$\{3,\ x,\ x^2,\ -3s+x^3,\ -4sx+x^4,\ -5sx^2+x^5,\ 3s^2-6sx^3+x^6,$$
$$7s^2x-7sx^4+x^7,\ 12s^2x^2-8sx^5+x^8,\ -3s^3+18s^2x^3-9sx^6+x^9\}$$

It is easy to guess that this sequence satisfies the recurrence

$$p_1(n+3,3,x,s) = xp_1(n+2,3,x,s) - sp_1(n,3,x,s) \qquad (2.9)$$

with initial values

$p_1(1,3,x,s)=x,\ p_1(2,3,x,s)=x^2,\ p_1(3,3,x,s)=x^3-3s.$

In order that the recurrence (2.9) remains true for $n=0$ we have set $p_1(0,3,x,s)=3.$

Further we get

$$\left(1-xz+sz^3\right)\sum_{n=0}^{\infty} p_1(n,3,x,s)z^n = \left(1-xz+sz^3\right)\sum_{n=0}^{3} p_1(n,3,x,s)z^n = 3-2xz.$$

This gives the conjecture

$$\sum_{n=0}^{\infty} p_1(n,3,x,s)z^n = \frac{3-2xz}{1-xz+sz^3}$$

which implies

$$\sum_{n\geq 0} p_1(n,3,x+1,x)z^n = \frac{3-2(x+1)z}{1-(x+1)z+xz^3} = \frac{1}{1-z}+\frac{2-xz}{1-xz-xz^2}.$$

By (2.7) the right-hand side is $\sum_{n\geq 0}\left(1+L_n(x,x)\right)z^n.$

Therefore we get

$$p_1(n,3,x+1,x)=1+L_n(x,x). \tag{2.10}$$

Since

$$\frac{1}{1-xz+sz^3}=\sum_{n\ge 0}z^n\left(x-sz^2\right)^n=\sum_{n,k}\binom{n}{k}z^{n+2k}x^{n-k}(-s)^k=\sum_n z^n\sum_{3k\le n}(-1)^k s^k\binom{n-2k}{k}x^{n-3k}$$

we get in analogy to (2.6) for $n>0$

$$p_1(n,3,x,s)=\sum_{3j\le n}(-1)^j\binom{n-2j}{j}\frac{n}{n-2j}x^{n-3j}s^j. \tag{2.11}$$

For

$$p_1(n,3,x,s)=\sum_{3k\le n}(-1)^k s^k\binom{n-2k}{k}x^{n-3k}-2s\sum_{3k\le n}(-1)^k s^k\binom{n-3-2k}{k}x^{n-3-3k}$$

$$=\sum_{3k\le n}(-1)^k s^k\left(\binom{n-2k}{k}+2\binom{n-1-2k}{k-1}\right)x^{n-3k}=\sum_{3k\le n}(-1)^k s^k\binom{n-2k}{k}\frac{n}{n-2k}x^{n-3k}.$$

In analogy to (2.4) and (2.8) we have

$$\sum_{k=0}^{\left\lfloor \frac{n}{3}\right\rfloor}\binom{n}{k}\tilde{p}_1(n-3k,3,x,s)s^k=x^n. \tag{2.12}$$

Observe that $\tilde{p}_1(n,3,x,s)+s\tilde{p}_1(n-3,3,x,s)=x\tilde{p}_1(n-1,3,x,s)$

We prove (2.12) by induction. It is true for $n=0$ and $n=1$. It is also true for $n=3$ because $\tilde{p}_1(3,3,x,s)+\binom{3}{1}s\tilde{p}_1(0,3,x,s)=x^3-3s+3s=x^3$. If it is true for $n-1\ge 3$ then

$$\sum_{k=0}^{\left\lfloor \frac{n}{3}\right\rfloor}\binom{n}{k}\tilde{p}_1(n-3k,3,x,s)s^k=\sum_k\binom{n-1}{k}\tilde{p}_1(n-3k,3,x,s)s^k+\sum_k\binom{n-1}{k-1}\tilde{p}_1(n-3k,3,x,s)s^k$$

$$=\sum_k\binom{n-1}{k}\left(\tilde{p}_1(n-3k,3,x,s)+s\tilde{p}_1(n-3k-3,3,x,s)\right)s^k=\sum_k\binom{n-1}{k}x\tilde{p}_1(n-1-3k,3,x,s)s^k=x^n.$$

The first terms of the sequence $\left(p_2(n,3,x,s)\right)_{n\ge 1}$ are

$\{0, -2sx, -3s^2, -2s^2x^2, -5s^3x, -3s^4-2s^3x^3, -7s^4x^2, -8s^5x-2s^4x^4, -3s^6-9s^5x^3\}$

Here we guess the recurrence relation

$$p_2(n+3,3,x,s)=sxzp_2(n+1,3,x,s)+s^2p_2(n,3,x,s)$$

which gives us the generating function

$$\sum_{n\ge 0}p_2(n,3,x,s)z^n=-\frac{3-xsz^2}{1-xsz^2-s^2z^3}. \tag{2.13}$$

Since

$$\sum_n (sz^2)^n (x+sz)^n = \sum_{n,k} (sz^2)^n \binom{n}{k} s^k z^k x^{n-k} = \sum_{n,k} s^{n+k} z^{2n+k} \binom{n}{k} x^{n-k} = \sum_n z^n \sum_j \binom{n-j}{2j-n} s^j x^{2n-3j}$$

we get

$$p_2(n,3,x,s) = -\sum_{3j\le 2n} \binom{n-j}{2n-3j} \frac{n}{n-j} s^j x^{2n-3j}. \tag{2.14}$$

This follows from

$$p_2(n,3,x,s) = xs\sum_j \binom{n-2-j}{2n-3j-4} s^j x^{2n-4-3j} - 3\sum_j \binom{n-j}{2n-3j} s^j x^{2n-3j}$$

$$= \sum_j \binom{n-1-j}{2n-3j-1} s^j x^{2n-3j} - 3\sum_j \binom{n-j}{2n-3j} s^j x^{2n-3j} = -\sum_j \binom{n-j}{2n-3j} \frac{n}{n-j} s^j x^{2n-3j}.$$

From (2.13) we conclude that

$$\sum_{n\ge 0} p_2(n,3,x+1,x) z^n = -\frac{3-x(x+1)z^2}{1-(x+1)xz^2-x^2z^3} = \frac{-1}{1+xz} - \frac{2-xz}{1-xz-xz^2}.$$

This gives

$$p_2(n,3,x+1,x) = (-1)^{n-1} x^n - L_n(x,x). \tag{2.15}$$

Comparing (2.10) and (2.15) we see that

$p_1(n,3,x+1,x) + p_2(n,3,x+1,x) = 1 + (-1)^{n-1} x^n,$

i.e. that (1.1) holds for $m=3$.

Therefore our guesses were correct and we have seen that

$$\sum_{3j\le n} (-1)^j \binom{n-2j}{j} \frac{n}{n-2j} (x+1)^{n-3j} x^j - \sum_{3j\le 2n} \binom{n-j}{2n-3j} \frac{n}{n-j} (x+1)^{2n-3j} x^j = 1 + (-1)^{n-1} x^n. \tag{2.16}$$

Before we consider the general case let us make some observations.

As we have seen we get some information about the coefficients $a(n,m,j)$ by obtaining recurrence relations for the polynomials $p_k(n,m,x,s)$.

If a sequence $a(n)$ satisfies the homogeneous linear recurrence $\sum_{k=0}^{m} b(n,k) a(n+k) = 0$ then we call $c(z) = \sum_{k=0}^{m} b(n,k) z^k$ a characteristic polynomial of the recurrence relation. Note that $c(z)$ is unique up to a multiplicative constant.

Let $f(z) = \sum_{k=0}^{d} v(k) z^k = (1-\alpha(1)z)(1-\alpha(2)z)\cdots(1-\alpha(d)z)$ with $\alpha(j) \ne 0$ for all $j$.

Then Newton's formula computes the power sums $pot(n) = \sum_{i=1}^{d} \alpha(i)^n$ in terms of the coefficients $v(k)$ as

$$\sum_{j=0}^{n-1} v(j)\, pot(n-j) + n v(n) = 0. \tag{2.17}$$

In particular we see that the sequence $\left(pot(n)\right)_{n\geq 1}$ satisfies the recurrence relation

$$\sum_{j=0}^{d} v(j)\,pot(n-j)=0 \tag{2.18}$$

for $n>m.$

Thus the characteristic polynomial of this recurrence is the reflected polynomial

$$c(z)=\sum_{j=0}^{d} v(d-j)z^{j}=\left(z-\alpha(1)\right)\cdots\left(z-\alpha(d)\right) \tag{2.19}$$

Let now $\overline{c}(z)$ be the characteristic polynomial of the recurrence of the power sums $\sum_{j=1}^{d}\left(\frac{b}{\alpha(j)}\right)^{n}$ for some constant $b.$

Then

$$\overline{c}(z)=\frac{(-z)^{d}}{v(d)}c\left(\frac{b}{z}\right). \tag{2.20}$$

This follows from

$$\left(z-\frac{b}{\alpha(1)}\right)\cdots\left(z-\frac{b}{\alpha(d)}\right)=\frac{1}{v(d)}\left(\alpha(1)z-b\right)\cdots\left(\alpha(d)z-b\right)$$
$$=\frac{(-z)^{d}}{v(d)}\left(\frac{b}{z}-\alpha(1)\right)\cdots\left(\frac{b}{z}-\alpha(d)\right).$$

## 3. Proof of Lemma 1

Using the concept of invariant polynomial I give a more direct approach to Lemma 1.

**a)** Let $f_{n,m}(x,s)=1-\prod_{j=0}^{n-1}\left(1-\omega^{j}x+\omega^{mj}s\right).$

We determine the part $p_k(n,m,x,s)$ of $f_{n,m}(x,s)$ with weight $k.$

As in [4] we write

$$1-xz+sz^{m}=\prod_{j=1}^{m}\left(1-u(j,m,x,s)z\right). \tag{3.1}$$

Then $1-\prod_{j=0}^{n-1}\left(1-\omega^{j}x+\omega^{mj}s\right)=1-\prod_{j=0}^{n-1}\prod_{k=1}^{m}\left(1-u(k,m,x,s)\omega^{j}\right)=\prod_{k=1}^{m}\prod_{j=0}^{n-1}\left(1-u(k,m,x,s)\omega^{j}\right).$

Therefore we get

$$1-\prod_{j=0}^{n-1}\left(1-\omega^{j}x+\omega^{mj}s\right)=1-\prod_{k=1}^{m}\prod_{j=0}^{n-1}\left(1-u(k,m,x,s)\omega^{j}\right)=1-\prod_{k=1}^{m}\left(1-u(k,m,x,s)^{n}\right).$$

This implies

$$f_{n,m}(x,s)=1-\prod_{j=0}^{n-1}\left(1-\omega^{j}x+\omega^{mj}s\right)=1-\prod_{k=1}^{m}\left(1-u(k,m,x,s)^{n}\right)=\sum_{k=1}^{m}p_k(n,m,x,s) \tag{3.2}$$

if we set

$$p_k(n,m,x,s)=(-1)^{k+1}\sum_{1\le i_1<i_2<\cdots<i_k\le m}\left(u(i_1,m,x,s)u(i_2,m,x,s)\cdots u(i_k,m,x,s)\right)^n. \qquad (3.3)$$

From

$$1-\lambda xz+\lambda^m sz^m=1-x(\lambda z)+s(\lambda z)^m=\prod_{j=1}^{m}\left(1-u(j,m,\lambda x,\lambda^m s)z\right)=\prod_{j=1}^{m}\left(1-u(j,m,x,s)\lambda z\right)$$

we see that by suitable ordering of the roots $u(j,x,s)$ we have

$$u(j,m,\lambda x,\lambda^m s)=\lambda u(j,m,x,s). \qquad (3.4)$$

Therefore $p_k(n,m,x,s)$ satisfies
$p_k(n,m,\lambda x,\lambda^m s)=\lambda^{kn}p_k(n,m,x,s).$
This means that $p_k(n,m,x,s)$ is the part of $f_{n,m}(x,s)$ with weight $k.$

**b)** The polynomial $f_{n,m}(x,s)$ has the form

$$f_{n,m}(x,s)=\sum_{k=1}^{m}p_k(n,m,x,s)=x^n+\sum_{k=1}^{m}\sum_{j=1+\left\lfloor\frac{(k-1)n}{m}\right\rfloor}^{\left\lfloor\frac{kn}{m}\right\rfloor}a(n,m,j)x^{kn-mj}s^j \qquad (3.5)$$

for some numbers $a(n,m,j).$

Let $x^is^j$ be of weight $k.$ This means that $i+mj=kn$ or $i=kn-mj.$ We must have $0\le j\le n$ and $0\le kn-mj\le n.$ Therefore $\frac{(k-1)n}{m}\le j\le\frac{kn}{m}.$
Thus we know that $p_k(n,m,x,s)$ is of the form
$p_k(n,m,x,s)=\sum_{\frac{(k-1)n}{m}\le j\le\frac{kn}{m}}a(n,m,j)x^{kn-mj}y^j.$ But (3.2) implies that $x^ns^j$ occurs in $f_{n,m}(x,s)$ if and only if $j=0.$
This gives

$$p_1(n,m,x,s)=x^n+\sum_{0<j\le\frac{n}{m}}a(n,m,j)x^{n-mj}y^j \qquad (3.6)$$

and

$$p_k(n,m,x,s)=\sum_{\frac{(k-1)n}{m}<j\le\frac{kn}{m}}a(n,m,j)x^{kn-mj}y^j \qquad (3.7)$$

for $k>1.$ Thus (3.5) is true.

(3.3) implies $p_n(n,m,x,s)=(-1)^{m+1}\left(u(1,m,x,s)u(2,m,x,s)\cdots u(m,m,x,s)\right)^n.$
Observing (3.1) this gives

$$p_n(n,m,x,s)=(-1)^{m+1+mn}s^n. \qquad (3.8)$$

Therefore (3.5) implies (1.1) if we change $x$ to $x+1$ and $s$ to $x.$

Consider for example $f_{6,4}(x,s)$.

Here we get $p_1(6,4,x,s)=x^6-6x^{6-4}s$, $p_2(6,4,x,s)=-3x^{12-8}s^2-2^{12-12}s^3$, $p_3(6,4,x,s)=3x^{18-16}s^4$, and $p_4(6,4,x,s)=-x^{24-24}s^6$.

Therefore $f_{6,4}(x,s)=x^6-6x^2s-3x^4s^2-2s^3+3x^2s^4-s^6$.

It is easily verified that $(x+1)^6-6(x+1)^2x-3(x+1)^4x^2-2x^3+3(x+1)^2x^4-x^6=1$.

**4) Some information about the coefficients** $a(n,m,j)$.

We first compute $p_1(n,m,x,s)=\sum_{j=1}^{m}u(j,m,x,s)^n$.

By (3.1) we have $f(z)=1-xz+sz^m$.

Newton's formula gives

$$p_1(n,m,x,s)-xp_1(n-1,m,x,s)+sp_1(n-m,m,x,s)=0 \qquad (4.1)$$

for $n>m$ and the initial values $p_1(n,m,x,s)=x^n$ for $n<m$ and $p_1(m,m,x,s)=x^m-ms$.

Thus a characteristic polynomial $c(m,1,x,s,z)$ of the recurrence of $p_1(n,m,x,s)$ is

$$c(m,1,x,s,z)=z^m-xz^{m-1}+s. \qquad (4.2)$$

As has been shown in [1] we have in analogy to (2.6) and (2.11)

$$p_1(n,m,x,s)=\sum_{j=0}^{\left\lfloor \frac{n}{m}\right\rfloor}(-1)^j\binom{n-(m-1)j}{j}\frac{n}{n-(m-1)j}x^{n-mj}s^j. \qquad (4.3)$$

This follows from (4.1). It is clear that the initial values coincide and the recurrence relation (4.1) is easily verified by comparing coefficients.

The same proof as for (2.12) gives

$$\sum_{k=0}^{\left\lfloor \frac{n}{m}\right\rfloor}\binom{n}{k}\tilde{p}_1(n-mk,m,x,s)s^k=x^n. \qquad (4.4)$$

In [1] it has also been shown that

$$(-1)^m p_{m-1}(n,m,x,(-1)^{m-1}s)=\sum_{j\le\frac{(m-1)n}{m}}\binom{n-j}{(m-1)n-mj}\frac{n}{n-j}x^{(m-1)n-mj}s^j \qquad (4.5)$$

with characteristic polynomial

$$c(m,m-1,x,s,z)=z^m-s^{m-2}xz+(-1)^m s^{m-1}. \qquad (4.6)$$

To prove this we see that by (3.3) the characteristic polynomial of the recurrence of the sequence $p_{m-1}(n,m,x,s)$ is

$$\bar{c}\left(m,1,x,s,z\right)=\frac{(-z)^m}{s}c\left(m,1,s,\frac{(-1)^m s}{z}\right)=\frac{(-z)^m}{s}\left(\left(\frac{(-1)^m s}{z}\right)^m-x\left(\frac{(-1)^m s}{z}\right)^{m-1}+s\right)$$

$$=(-1)^m\left((-1)^m s^{m-1}-s^{m-2}xz+z^m\right).$$

Furthermore we get

$p_{m-1}(n,m,x,s)=0$ for $n<m-1$, $p_{m-1}(m-1,m,x,s)=(-1)^m(m-1)s^{m-2}x$, $p_{m-1}(m,m,x,s)=-ms^{m-1}$.

For $n>m$ (4.6) implies

$$p(n,m,x,s)-s^{m-2}xzp(n-m+1,m,x,s)-(-s)^{m-1}p(n-m,m,x,s)=0.$$

Comparing coefficients we get (4.5).

More generally formula (2.20) implies

$$c(m,m-k,x,s,z)=\frac{z^{\binom{m}{k}}}{s^{\binom{m-1}{k-1}}}c\left(m,k,x,s,(-1)^m\frac{s}{z}\right). \qquad (4.7)$$

For in this case

$$f(z)=\sum_{k=0}^{d}v(k)z^k=\prod_{i_1<i_2<\cdots<i_k}\left(1-u_{i_1}(m,x,s)u_{i_2}(m,x,s)\cdots u_{i_k}(m,x,s)z\right).$$

Thus $d=\binom{m}{k}$, $v(d)=(-1)^{\binom{m}{k}}s^{\binom{m-1}{k-1}}$ and $b=(-1)^m s.$

In particular we see that the sequences $\left(p_k(n,m,x,s)\right)$ satisfy a linear recurrence of order $\binom{m}{k}$ and of no smaller order.

Except for $k=1$ and $k=m-1$ we have neither for the polynomials $p_k(n,m,x,s)$ nor for their characteristic polynomials $c(m,k,x,s,z)$ explicit formulae which are valid for all $m$. For $2\le m\le 8$ we have computed all characteristic polynomials but found no conjecture for their general form.

By substituting $(x,s)\to(x+1,x)$ we get

$$\prod_{j=1}^{m}\left(1-u(j,m,x+1,x)z\right)=1-(x+1)z+xz^m=\left(1-z\right)\left(1-xz\frac{1-z^{m-1}}{1-z}\right).$$

In this case one root $u(1,m,x+1,x)$ is 1.

Therefore the characteristic polynomial $c(m,k,x,s,z)$ splits into the product of two polynomials with integer coefficients $c(m,k,x+1,x,z)=w_{k-1}(m,x,z)w_k(m,x,z).$

The first one has roots $u(1,m,x+1,x)u(i_1,m,x+1,x)\cdots u\left(i_{k-1},m,x+1,x\right)$ and the other all products $u\left(i_1,m,x+1,x\right)\cdots u\left(i_k,m,x+1,x\right)$ with $1<i_1<\cdots<i_k$. This observation proves a conjecture in [2].

For small values of $m$ it is easy to find the recurrence relations for $p_k(n,m,x,s)$.
Since we know that $\left(p_k(n,m,x,s)\right)_{n\ge 0}$ satisfies a linear recurrence $c(m,k,x,s,z)$ of order $\binom{m}{k}$ it suffices to guess such a recurrence and verify it for small values of $n$.

For example we guess that $c(4,2,x,s,z)=z^6-sz^4-sx^2z^3-s^2z^2+s^3$.
Since $z^6c\left(4,2,x,s,\frac{1}{z}\right)=1-sz^2-sx^2z^3-s^2z^4+s^3z^6$ and

$$(1-sz^2-sx^2z^3-s^2z^4+s^3z^6)\sum_{n=0}^{15}p_2(n,4,x,s)z^n=-2sz^2-3x^2sz^3-4s^2z^4+6s^3z^6+O(z^{16})$$

we see that our guess is correct.

It should be noted that $c(m,k,x,s,z)$ contains the whole information about the polynomials $p_k(n,m,x,s)$, because the polynomials $(-1)^{k-1}p_k(n,m,x,s)$ are the power sums of the reflected polynomial $c^*(m,k,x,s,z)=z^{\binom{m}{k}}c\left(m,k,x,s,\frac{1}{z}\right)$.
By (3.4) we get

$$c^*(m,k,ax,a^ms,z)=c^*\left(m,k,x,s,a^kz\right). \tag{4.8}$$

This means that $c^*\left(m,k,x,s,z\right)=\sum_{0\le n\le\binom{m}{k},0\le j\le\frac{nk}{m}}h(j)x^{nk-mj}s^jz^n$

for some coefficients $h(j)$.

## 5. The original problem revisited

Consider now the linear operators $E$ and $\Delta$ on the vector space of polynomials defined by $Ef(x)=f(x+1)$ and $\Delta f(x)=(E-I)f(x)=f(x+1)-f(x)$, where $I$ is the identity $If(x)=f(x)$.
Then $f_{i,m}(E,\Delta)=I$ or equivalently

$$E^i+\sum_{k=1}^{m-1}\sum_{j=1+\left\lfloor\frac{(k-1)i}{m}\right\rfloor}^{\left\lfloor\frac{ki}{m}\right\rfloor}a(i,m,j)E^{ki-mj}\Delta^j=I+(-1)^{m(i-1)}\Delta^i. \tag{5.1}$$

For the polynomials $\binom{x}{r}$ we have $E^k\binom{x}{r}=\binom{x+k}{r}$ and $\Delta^k\binom{x}{r}=\binom{x}{r-k}$.

If we apply the operator (5.1) to $\binom{x}{r}$ we get

$$\binom{x+i}{r}+\sum_{k=1}^{m-1}\sum_{j=1+\left\lfloor\frac{(k-1)i}{m}\right\rfloor}^{\left\lfloor\frac{ki}{m}\right\rfloor}a(i,m,j)\binom{x+ki-mj}{r-j}=\binom{x}{r}+(-1)^{m(i-1)}\binom{x}{r-i}.$$

Let now $r=\left\lfloor\frac{n+ih+\ell)}{m}\right\rfloor$.
Then we get

$$\binom{n}{\left\lfloor \frac{n+ih+\ell)}{m} \right\rfloor}+(-1)^{m(i-1)}\binom{n}{\left\lfloor \frac{n+i(h-m)+\ell)}{m} \right\rfloor}$$

$$=\binom{n+i}{\left\lfloor \frac{n+ih+\ell)}{m} \right\rfloor}+\sum_{k=1}^{m-1}\sum_{j=1+\left\lfloor \frac{(k-1)i}{m}\right\rfloor}^{\left\lfloor \frac{ki}{m}\right\rfloor} a(i,m,j)\binom{n+ki-mj}{\left\lfloor \frac{n+ih+\ell-jm)}{m} \right\rfloor}.$$

We multiply this identity by $z^h$ and sum over all $h\in\mathbb{Z}$ and set

$$A(n,m,i,\ell,z)=\sum_{h\in\mathbb{Z}}\binom{n}{\left\lfloor \frac{n+ih+\ell)}{m} \right\rfloor}z^h.$$

Note that this is a finite sum.
Since

$$\sum_{h\in\mathbb{Z}}\binom{n}{\left\lfloor \frac{n+i(h-m)+\ell)}{m} \right\rfloor}z^h = z^m\sum_{h\in\mathbb{Z}}\binom{n}{\left\lfloor \frac{n+i(h-m)+\ell)}{m} \right\rfloor}z^{h-m} = z^m A(n,m,i,\ell,z)$$

$$\sum_{h\in\mathbb{Z}}\binom{n+ki-mj}{\left\lfloor \frac{n+ih-jm+\ell)}{m} \right\rfloor}z^h = \sum_{h\in\mathbb{Z}}\binom{n+ki-mj}{\left\lfloor \frac{n+ki-mj+i(h-k)+\ell)}{m} \right\rfloor}z^h = z^k A(n+ki-mj,m,i,\ell,z)$$

we get

$$A(n,m,i,\ell,z)+(-1)^{i(m-1)}z^m A(n,m,i,\ell,z)=zA(n+i,m,i,\ell,z)$$

$$+\sum_{k=1}^{m-1}\sum_{j=1+\left\lfloor \frac{(k-1)i}{m}\right\rfloor}^{\left\lfloor \frac{ki}{m}\right\rfloor} a(i,m,j)z^k A(n+ki-mj,m,i,\ell,z).$$

Thus we have the main result of [1] and [2].

**Theorem 1**

*Let* $m\geq 2, i\geq 1$ *be integers,* $n\in\mathbb{N}$ *and* $\ell\in\mathbb{Z}$.
*The sequences*

$$A(n,m,i,\ell,z)=\sum_{h\in\mathbb{Z}}\binom{n}{\left\lfloor \frac{n+ih+\ell}{m} \right\rfloor}z^h \in\mathbb{Q}\left[z,z^{-1}\right] \tag{5.2}$$

*satisfy the linear recurrence of order* $i$ *with constant integer coefficients*

$$\sum_{k=1}^{m-1} z^{k-1}p_k(i,m,E,1)A(n,m,i,\ell,z)=\left(\frac{1}{z}+(-1)^{m(i-1)}z^{m-1}\right)A(n,m,i,\ell,z). \tag{5.3}$$

As an example consider the case $m=2$, $i=5$ and $\ell=0$.

The first terms of $A(n,2,5,0,z)$ are
$\left\{1, 1, 2, 3, 6+z, 10+\frac{1}{z}+z, 20+\frac{1}{z}+6z, 35+\frac{7}{z}+7z, 70+\frac{8}{z}+28z, 126+\frac{36}{z}+36z+z^2\right\}$
We have $p_1(5,2,x,1)=x^5-5x^3+5x$.
Therefore $A(n,2,5,0,z)$ satisfies the recurrence

$$A(n+5,2,5,0,z)-5A(n+3,2,5,0,z)+5A(n+1,2,5,0,z)-\left(z+\frac{1}{z}\right)A(n,2,5,0,z)=0.$$

For $z=-1$ the sequence begins with
$\{1, 1, 2, 3, 5, 8, 13, 21, 34\}$
and the recurrence reduces to

$$A(n+5,2,5,0,-1)-5A(n+3,2,5,0,-1)+5A(n+1,2,5,0,-1)+2A(n,2,5,0,-1)=0.$$

The characteristic polynomial of this recurrence is
$z^5-5z^3+5z+2=(z+2)\left(z^2-z-1\right)^2$.
Since the characteristic polynomials of the recurrence of the Fibonacci numbers $F_{n+1}$ is $z^2-z-1$ and the first 5 terms of the sequence $A(n,2,5,0,-1)$ coincide with $F_{n+1}$ we have
$A(n,2,5,0,-1)=F_{n+1}$.
With other words we have

$$F_{n+1}=\sum_{h\in\mathbb{Z}}(-1)^h\binom{n}{\left\lfloor \frac{n+5h}{2}\right\rfloor}. \qquad (5.4)$$

In the same way we get

$$F_{n}=\sum_{h\in\mathbb{Z}}(-1)^h\binom{n}{\left\lfloor \frac{n+5h+2}{2}\right\rfloor}. \qquad (5.5)$$

These two curious identities which are due to Issai Schur [5] in a more general form were the starting point of my 2001 paper [1].